\documentclass[onecolumn]{autart}
\usepackage{amsmath}
\begin{document}
\input{amssym}
\begin{frontmatter}
\title{Group classification of steady two-dimensional boundary-layer stagnation-point flow equations}
\author[MN]{Mehdi Nadjafikhah}\ead{m\_nadjafikhah@iust.ac.ir},
\author[SRH]{Seyed Reza Hejazi}\ead{reza\_hejazi@iust.ac.ir}
\address[MN]{School of Mathematics, Iran University of Science and Technology, Narmak-16, Teharn, I.R.Iran.}
\address[SRH]{Same as author one.}
\begin{abstract}
Lie symmetry group method is applied to study the boundary-layer
equations for two-dimensional steady flow of an incompressible,
viscous fluid near a stagnation point at a heated stretching sheet
placed in a porous medium equation. The symmetry group and its
optimal system are given, and group invariant solutions associated
to the symmetries are obtained. Finally the structure of the Lie
algebra symmetries is determined.
\end{abstract}
\end{frontmatter}
\section{Introduction}
In fluid dynamics, a stagnation point is a point in a flow field
where the local velocity of the fluid is zero. Stagnation points
exist at the surface of objects in the flow field, where the fluid
is brought to rest by the object. The Bernoulli equation shows
that the static pressure is highest when the velocity is zero and
hence static pressure is at its maximum value at stagnation
points. This static pressure is called the stagnation pressure.
The Bernoulli equation applicable to incompressible flow shows
that the stagnation pressure is equal to the dynamic pressure plus
static pressure. Total pressure is also equal to dynamic pressure
plus static pressure so, in incompressible flows, stagnation
pressure is equal to total pressure. (In compressible flows,
stagnation pressure is also equal to total pressure providing the
fluid entering the stagnation point is brought to rest is
entropically.)

In physics and fluid mechanics, a boundary layer is that layer of
fluid in the immediate vicinity of a bounding surface. In the
Earth's atmosphere, the planetary boundary layer is the air layer
near the ground affected by diurnal heat, moisture or momentum
transfer to or from the surface. On an aircraft wing the boundary
layer is the part of the flow close to the wing. The boundary
layer effect occurs at the field region in which all changes occur
in the flow pattern. The boundary layer distorts surrounding
non-viscous flow. It is a phenomenon of viscous forces. This
effect is related to the Reynolds number (In fluid mechanics and
heat transfer, the Reynold's number is a dimensionless number that
gives a measure of the ratio of inertial forces to viscous and,
consequently, it quantifies the relative importance of these two
types of forces for given flow conditions). Laminar boundary
layers come in various forms and can be loosely classified
according to their structure and the circumstances under which
they are created. The thin shear layer which develops on an
oscillating body is an example of a Stokes boundary layer, whilst
the Blasius boundary layer refers to the well-known similarity
solution for the steady boundary layer attached to a flat plate
held in an oncoming unidirectional flow. When a fluid rotates,
viscous forces may be balanced by the Coriolis effect, rather than
convective inertia, leading to the formation of an Ekman layer.
Thermal boundary layers also exist in heat transfer. Multiple
types of boundary layers can coexist near a surface
simultaneously. The deduction of the boundary layer equations was
perhaps one of the most important advances in fluid dynamics.
Using an order of magnitude analysis, the well-known governing
Navier됩tokes equations of viscous fluid flow can be greatly
simplified within the boundary layer. Notably, the characteristic
of the partial differential equations (PDE) becomes parabolic,
rather than the elliptical form of the full Navier됩tokes
equations. This greatly simplifies the solution of the equations.
By making the boundary layer approximation, the flow is divided
into an inviscid portion (which is easy to solve by a number of
methods) and the boundary layer, which is governed by an easier to
solve PDE.

Flow and heat transfer of an incompressible viscous fluid over a
stretching sheet appear in several manufacturing processes of
industry such as the extrusion of polymers, the cooling of
metallic plates, the aerodynamic extrusion of plastic sheets, etc.
In the glass industry, blowing, floating or spinning of fibres are
processes, which involve the flow due to a stretching surface.
Mahapatra and Gupta studied the steady two-dimensional
stagnation-point flow of an incompressible viscous fluid over a
flat deformable sheet when the sheet is stretched in its own plane
with a velocity proportional to the distance from the
stagnation-point. They concluded that, for a fluid of small
kinematic viscosity, a boundary layer is formed when the
stretching velocity is less than the free stream velocity and an
inverted boundary layer is formed when the stretching velocity
exceeds the free stream velocity. Temperature distribution in the
boundary layer is determined when the surface is held at constant
temperature giving the so called surface heat flux. In their
analysis, they used the finite-differences scheme along with the
Thomas algorithm to solve the resulting system of ordinary
differential equations.

This paper is concerned with the solution of steady
two-dimensional stagnation point flow of an incompressible viscous
fluid over a stretching sheet which is placed in a fluid saturated
porous medium. Lie-group theory is applied to the equations of
motion for determining symmetry reductions of partial differential
equations.
\section{Lie Symmetries of the Equations}
A PDE with $p-$independent and $q-$dependent variables has a Lie
point transformations
\begin{eqnarray*}
\widetilde{x}_i=x_i+\varepsilon\xi_i(x,u)+{\mathcal
O}(\varepsilon^2),\qquad
\widetilde{u}_{\alpha}=u_\alpha+\varepsilon\varphi_\alpha(x,u)+{\mathcal
O}(\varepsilon^2),
\end{eqnarray*}
where
$\displaystyle{\xi_i=\frac{\partial\widetilde{x}_i}{\partial\varepsilon}\Big|_{\varepsilon=0}}$
for $i=1,...,p$ and
$\displaystyle{\varphi_\alpha=\frac{\partial\widetilde{u}_\alpha}{\partial\varepsilon}\Big|_{\varepsilon=0}}$
for $\alpha=1,...,q$. The action of the Lie group can be
considered by its associated infinitesimal generator
\begin{eqnarray}\label{eq:18}
\textbf{v}=\sum_{i=1}^p\xi_i(x,u)\frac{\partial}{\partial{x_i}}+\sum_{\alpha=1}^q\varphi_\alpha(x,u)\frac{\partial}{\partial{u_\alpha}}
\end{eqnarray}
on the total space of PDE (the space containing independent and
dependent variables). Furthermore, the characteristic of the
vector field (\ref{eq:18}) is given by
\begin{eqnarray*}
Q^\alpha(x,u^{(1)})=\varphi_\alpha(x,u)-\sum_{i=1}^p\xi_i(x,u)\frac{\partial
u^\alpha}{\partial x_i},
\end{eqnarray*}
and its $n-$th prolongation is determined by
\begin{eqnarray*}
\textbf{v}^{(n)}=\sum_{i=1}^p\xi_i(x,u)\frac{\partial}{\partial
x_i}+\sum_{\alpha=1}^q\sum_{\sharp
J=j=0}^n\varphi^J_\alpha(x,u^{(j)})\frac{\partial}{\partial
u^\alpha_J},
\end{eqnarray*}
where
$\varphi^J_\alpha=D_JQ^\alpha+\sum_{i=1}^p\xi_iu^\alpha_{J,i}$.
($D_J$ is the total derivative operator describes in
(\ref{eq:19})).

The deduction of the boundary layer equations was perhaps one of
the most important advances in fluid dynamics. Using an order of
magnitude analysis, the well-known governing Navier됩tokes
equations of viscous fluid flow can be greatly simplified within
the boundary layer. Notably, the characteristic of the partial
differential equations becomes parabolic, rather than the
elliptical form of the full Navier됩tokes equations. This greatly
simplifies the solution of the equations. By making the boundary
layer approximation, the flow is divided into an inviscid portion
(which is easy to solve by a number of methods) and the boundary
layer, which is governed by an easier to solve PDE. The continuity
and Navier됩tokes equations for a two-dimensional steady
incompressible flow in Cartesian coordinates are given by
\begin{eqnarray}\label{eq:1} \left\{ \begin{array}{l}
\displaystyle \frac{\partial U}{\partial x}+\frac{\partial V}{\partial
y}=0,\\[3mm]
\displaystyle U\frac{\partial U}{\partial x}+V\frac{\partial
V}{\partial y}=P\frac{\partial P}{\partial x}+\nu\frac{\partial^2
U}{\partial y^2}+\frac{\nu}{k}(P-U),\\[3mm]
\displaystyle U\frac{\partial T}{\partial x}+V\frac{\partial
T}{\partial y}=\alpha\frac{\partial^2 T}{\partial y^2}.\end{array}
\right.
\end{eqnarray}
where $U,V,P$ and $T$ are smooth functions of
$\displaystyle{(x,y)}$, and also $U$ and $V$ are the velocity
components along $x$-axis and $y$-axis respectively, $\nu$ is the
kinematic viscosity, $k$ is the permeability of the porous medium,
$T$ is the fluid temperature, $\alpha$ is the coefficient of
thermal diffusivity. The aim is to analysis the Lie point symmetry
structure of the system of steady two-dimensional boundary-layer
stagnation-point flow equations.

Let us consider a one-parameter Lie group of infinitesimal
transformations $(x,y,U,V,P,T)$ given by
\begin{align*}
\widetilde{x}&=x+\varepsilon\xi_1(x,y,U,V,P,T)+{\mathcal
O}(\varepsilon^2),&
&&\widetilde{y}&=y+\varepsilon\xi_2(x,y,U,V,P,T)+{\mathcal
O}(\varepsilon^2),\\
\widetilde{U}&=U+\varepsilon\eta_1(x,y,U,V,P,T)+{\mathcal
O}(\varepsilon^2),&
&&\widetilde{V}&=V+\varepsilon\eta_2(x,y,U,V,P,T)+{\mathcal
O}(\varepsilon^2),\\
\widetilde{P}&=P+\varepsilon\eta_3(x,y,U,V,P,T)+{\mathcal
O}(\varepsilon^2),&
&&\widetilde{T}&=T+\varepsilon\eta_4(x,y,U,V,P,T)+{\mathcal
O}(\varepsilon^2),
\end{align*}
where $\varepsilon$ is the group parameter. Then one requires that
this transformations leaves invariant the set of solutions of the
system (\ref{eq:1}). This yields to the linear system of equations
for the infinitesimals $\xi_1(x,y,U,V,P,T),\xi_2(x,y,U,V,P,T)$,
$\eta_1(x,y,U,V,P,T),\eta_2(x,y,U,V,P,T),\eta_3(x,y,U,V,P,T)$ and
$\eta_4(x,y,U,V,P,T)$. The Lie algebra of infinitesimal symmetries
is the set of vector fields in the form of
\begin{eqnarray*}
\textbf{v}=\xi_1\frac{\partial}{\partial
x}+\xi_2\frac{\partial}{\partial y}+\eta_1\frac{\partial}{\partial
U}+\eta_2\frac{\partial}{\partial
V}+\eta_3\frac{\partial}{\partial
P}+\eta_4\frac{\partial}{\partial T}.
\end{eqnarray*}
This vector field has the second prolongation
\begin{eqnarray}\label{eq:2}
\textbf{v}^{(2)}&=&\textbf{v}+\varphi_1^x\frac{\partial}{\partial
x}+\varphi_1^y\frac{\partial}{\partial
y}+\varphi_1^{xx}\frac{\partial}{\partial
U_{xx}}+\varphi_1^{xy}\frac{\partial}{\partial
U_{xy}}+\varphi_1^{yy}\frac{\partial}{\partial U_{yy}}
+\varphi_2^x\frac{\partial}{\partial
V_x}+\varphi_2^y\frac{\partial}{\partial
V_y}+\varphi_2^{xx}\frac{\partial}{\partial
V_{xx}}\nonumber\\
&&+\varphi_2^{xy}\frac{\partial}{\partial
V_{xy}}+\varphi_2^{yy}\frac{\partial}{\partial
V_{yy}}+\varphi_3^x\frac{\partial}{\partial
P_x}+\varphi_3^y\frac{\partial}{\partial
P_y}+\varphi_3^{xx}\frac{\partial}{\partial
P_{xx}}+\varphi_3^{xy}\frac{\partial}{\partial
P_{xy}}+\varphi_3^{yy}\frac{\partial}{\partial
P_{yy}}\\
&&+\varphi_4^x\frac{\partial}{\partial
T_x}+\varphi_4^y\frac{\partial}{\partial
T_y}+\varphi_4^{xx}\frac{\partial}{\partial T_{xx}}
+\varphi_4^{xy}\frac{\partial}{\partial
T_{xy}}+\varphi_4^{yy}\frac{\partial}{\partial T_{yy}},\nonumber
\end{eqnarray}
with the coefficients
\begin{align*}
\varphi_1^x&=D_x(\varphi_1-\xi_1U_x-\xi_2U_y)+\xi_1U_{xx}+\xi_2U_{xy},&
\varphi_1^y&=D_y(\varphi_1-\xi_1U_x-\xi_2U_y)+\xi_1U_{xy}+\xi_2U_{yy},\\
\varphi_1^{xx}&=D_x^2(\varphi_1-\xi_1U_x-\xi_2U_y)+\xi_1U_{xxx}+\xi_2U_{xxy},&
\varphi_1^{xy}&=D_xD_y(\varphi_1-\xi_1U_x-\xi_2U_y)+\xi_1U_{xxy}+\xi_2U_{xxy},\\
\varphi_1^{yy}&=D_y^2(\varphi_1-\xi_1U_x-\xi_2U_y)+\xi_1U_{xxy}+\xi_2U_{xyy},&
\varphi_2^x&=D_x(\varphi_2-\xi_1V_x-\xi_2V_y)+\xi_1V_{xx}+\xi_2V_{xy},\\
\varphi_2^y&=D_y(\varphi_2-\xi_1V_x-\xi_2V_y)+\xi_1V_{xy}+\xi_2V_{yy},&
\varphi_2^{xx}&=D_x^2(\varphi_2-\xi_1V_x-\xi_2V_y)+\xi_1V_{xxx}+\xi_2V_{xxy},\\
\varphi_2^{xy}&=D_xD_y(\varphi_2-\xi_1V_x-\xi_2V_y)+\xi_1V_{xxy}+\xi_2V_{xxy},&
\varphi_2^{yy}&=D_y^2(\varphi_2-\xi_1P_x-\xi_2P_y)+\xi_1P_{xxy}+\xi_2P_{xyy},\\
\varphi_3^x&=D_x(\varphi_3-\xi_1P_x-\xi_2P_y)+\xi_1P_{xx}+\xi_2P_{xy},&
\varphi_3^y&=D_y(\varphi_3-\xi_1P_x-\xi_2P_y)+\xi_1P_{xy}+\xi_2P_{yy},\\
\varphi_3^{xx}&=D_x^2(\varphi_3-\xi_1P_x-\xi_2P_y)+\xi_1P_{xxx}+\xi_2P_{xxy},&
\varphi_3^{xy}&=D_xD_y(\varphi_3-\xi_1T_x-\xi_2T_y)+\xi_1T_{xxy}+\xi_2T_{xxy},\\
\varphi_3^{yy}&=D_y^2(\varphi_3-\xi_1T_x-\xi_2T_y)+\xi_1T_{xxy}+\xi_2T_{xyy},&
\varphi_4^x&=D_x(\varphi_4-\xi_1T_x-\xi_2T_y)+\xi_1T_{xx}+\xi_2T_{xy},\\
\varphi_4^y&=D_y(\varphi_4-\xi_1T_x-\xi_2T_y)+\xi_1T_{xy}+\xi_2T_{yy},&
\varphi_4^{xx}&=D_x^2(\varphi_4-\xi_1P_x-\xi_2P_y)+\xi_1P_{xxx}+\xi_2P_{xxy},\\
\varphi_4^{xy}&=D_xD_y(\varphi_4-\xi_1T_x-\xi_2T_y)+\xi_1T_{xxy}+\xi_2T_{xxy},&
\varphi_4^{yy}&=D_y^2(\varphi_4-\xi_1T_x-\xi_2T_y)+\xi_1T_{xxy}+\xi_2T_{xyy},
\end{align*}
where the operators $D_x$ and $D_y$ denote the total derivative
with respect to $x$ and $t$:
\begin{eqnarray}\label{eq:19}
D_x=\frac{\partial}{\partial x}+U_x\frac{\partial}{\partial
U}+V_x\frac{\partial}{\partial V}+P_x\frac{\partial}{\partial
P}+T_x\frac{\partial}{\partial T}\cdots,\quad
D_y=\frac{\partial}{\partial y}+U_y\frac{\partial}{\partial
U}+V_y\frac{\partial}{\partial V}+P_y\frac{\partial}{\partial
P}+T_y\frac{\partial}{\partial T}\cdots.
\end{eqnarray}
Using the invariance condition, i.e., applying the second
prolongation $\textbf{v}^{(2)}$, (\ref{eq:2}), to system
(\ref{eq:1}), and by solving the linear system
\begin{eqnarray*} \left\{\begin{array}{l} \displaystyle
\textbf{v}^{(2)}\Big(\frac{\partial U}{\partial x}+\frac{\partial
V}{\partial y}\Big)=0,\\[3mm] \displaystyle \textbf{v}^{(2)}\Big(U\frac{\partial
U}{\partial x}+V\frac{\partial V}{\partial y}-P\frac{\partial
P}{\partial x}-\nu\frac{\partial^2 U}{\partial
y^2}-\frac{\nu}{k}(P-U)\Big)=0,\\[3mm]
\displaystyle \textbf{v}^{(2)}\Big(U\frac{\partial T}{\partial
x}+V\frac{\partial T}{\partial y}-\alpha\frac{\partial^2
T}{\partial y^2}\Big)=0,\end{array}\right.\quad
\pmod{(\ref{eq:1})}
\end{eqnarray*}
the following system of 112 determining equations yields:
\begin{eqnarray*}
&\displaystyle 2\alpha\frac{\partial^2\xi_2}{\partial T\partial
y}+2V\frac{\partial\xi_2}{\partial
T}-\alpha\frac{\partial^2\eta_4}{\partial T^2}=0,\qquad \alpha
V\frac{\partial\xi_1}{\partial
U}+\alpha\nu\frac{\partial^2\xi_1}{\partial V\partial y}-\nu
V\frac{\partial\xi_1}{\partial V}+2\nu
U\frac{\partial\xi_2}{\partial V}-\alpha
U\frac{\partial\xi_1}{\partial U}+\nu
U\frac{\partial\xi_1}{\partial U}=0,\\
&\displaystyle V\frac{\partial\xi_1}{\partial
U}+U\frac{\partial\xi_2}{\partial V}=0,\qquad
U\frac{\partial\xi_1}{\partial V}+V\frac{\partial\xi_1}{\partial
V}=0,\quad \dots\dots\quad 2\frac{\partial^2\xi_2}{\partial
U\partial y}-\frac{\partial^2\eta_1}{\partial U^2}=0,\qquad
\frac{\partial^2\xi_1}{\partial
U^2}-2\frac{\partial\xi_2}{\partial^2 V\partial U}=0.
\end{eqnarray*}
The solution of the above system gives the following coefficients
of the vector field $\textbf{v}$:
\begin{eqnarray*}
\begin{array}{lclclclclcl}
\xi_1=C_1,&&\xi_2=C_2,&&\eta_1=0,&&\eta_2=0,&&\eta_3=T,&&\eta_4=0,
\end{array}
\end{eqnarray*}
where $C_1$ and $C_2$ are arbitrary constants, thus the Lie
algebra ${\goth g}$ of the steady two-dimensional boundary-layer
stagnation-point flow equations is spanned by the four vector
fields
\begin{eqnarray*}
\begin{array}{lclclcl}
\displaystyle{\textbf{v}_1=\frac{\partial}{\partial x}},&&
\displaystyle{\textbf{v}_2=\frac{\partial}{\partial y}},&&
\displaystyle{\textbf{v}_3=\frac{\partial}{\partial
T}},&&\displaystyle{\textbf{v}_4=T\frac{\partial}{\partial T}},
\end{array}
\end{eqnarray*}
which $\textbf{v}_1,\textbf{v}_2$ and $\textbf{v}_3$ are
translation on $x,u$ and $T$, $\textbf{v}_4$ is scaling on $T$.
The commutation relations between these vector fields is given by
the table 1, where entry in row $i$ and column $j$ representing
$[\textbf{v}_i,\textbf{v}_j]$.
\begin{table}
\caption{Commutation relations of $\goth g$ }\label{table:1}
\vspace{-0.3cm}\begin{eqnarray*}\hspace{-0.75cm}\begin{array}{lclclclclc}
\hline
  [\,,\,]       &\hspace{1.1cm}\textbf{v}_1  &\hspace{0.5cm}\textbf{v}_2  &\hspace{0.5cm}\textbf{v}_3          &\hspace{0.5cm}\textbf{v}_4  \\ \hline
  \textbf{v}_1  &\hspace{1.1cm} 0            &\hspace{0.5cm} 0            &\hspace{0.5cm}0                     &\hspace{0.5cm}0  \\
  \textbf{v}_2  &\hspace{1.1cm} 0            &\hspace{0.5cm} 0            &\hspace{0.5cm}0                     &\hspace{0.5cm}0  \\
  \textbf{v}_3  &\hspace{1.1cm} 0            &\hspace{0.5cm} 0            &\hspace{0.5cm}0                     &\hspace{0.5cm}\textbf{v}_3 \\
  \textbf{v}_4  &\hspace{1.1cm} 0            &\hspace{0.5cm} 0            &\hspace{0.5cm}-\textbf{v}_3         &\hspace{0.5cm}0  \\
  \hline\end{array}\end{eqnarray*}\end{table}

The one-parameter groups $G_i$ generated by the base of $\goth g$
are given in the following table.
\begin{eqnarray*}
\begin{array}{lclclclc}
G_1:(x+\varepsilon,y,U,V,P,T),&& G_2:(x,y+\varepsilon,U,V,P,T),&&
G_3:(x,y,U,V,P,T+\varepsilon),&& G_4:(x,y,U,V,P,Te^\varepsilon).
\end{array}
\end{eqnarray*}
Since each group $G_i$ is a symmetry group and if
$U=f(x,y),V=g(x,y),P=h(x,y)$ and $T=r(x,y)$ are solutions of the
system (\ref{eq:1}), so are the functions
\begin{align*}
U_1&=f(x+\varepsilon,y),&U_1&=f(x,y+\varepsilon),&U_1&=f(x,y),&U_1&=f(x,y),\\
V_1&=g(x+\varepsilon,y),&V_1&=g(x,y+\varepsilon),&V_1&=g(x,y),&V_1&=g(x,y),\\
P_1&=h(x+\varepsilon,y),&P_1&=h(x,y+\varepsilon),&P_1&=h(x,y),&P_1&=h(x,y),\\
T_1&=r(x+\varepsilon,y),&T_1&=r(x,y+\varepsilon),&T_1&=r(x,y)+\varepsilon,&T_1&=e^{-\varepsilon}r(x,y),
\end{align*}
where $\varepsilon$ is a real number. Here we can find the general
group of the symmetries by considering a general linear
combination $c_1\textbf{v}_1+\cdots+c_1\textbf{v}_4$ of the given
vector fields. In particular if $g$ is the action of the symmetry
group near the identity, it can be represented in the form
$g=\exp(\varepsilon_4\textbf{v}_4)\circ\cdots\circ\exp(\varepsilon_1\textbf{v}_1)$.
\section{Symmetry reduction for steady two-dimensional boundary-layer stagnation-point flow equations}
Lie-group method is applicable to both linear and non-linear
partial differential equations, which leads to similarity
variables that may be used to reduce the number of independent
variables in partial differential equations. By determining the
transformation group under which a given partial differential
equation is invariant, we can obtain information about the
invariants and symmetries of that equation.

The first advantage of symmetry group method is to construct new
solutions from known solutions. Neither the first advantage nor
the second will be investigated here, but symmetry group method
will be applied to the Eq. (\ref{eq:1}) to be connected directly
to some order differential equations. To do this, a particular
linear combinations of infinitesimals are considered and their
corresponding invariants are determined. The steady
two-dimensional boundary-layer stagnation-point flow equations
expressed in the coordinates $(x,y)$, so to reduce this equation
is to search for its form in specific coordinates. Those
coordinates will be constructed by searching for independent
invariants $(r,s)$ corresponding to an infinitesimal generator. So
using the chain rule, the expression of the equation in the new
coordinate allows us to the reduced equation. Here we will obtain
some invariant solutions with respect to symmetries. First we
obtain the similarity variables for each term of the Lie algebra
$\goth g$, then we use this method to reduced the PDE and find the
invariant solutions. Here our system has two-independent and
four-dependent variables, thus, the similarity transformations and
invariant functions with invariant solutions are coming in the
following table.
\begin{eqnarray*}
\begin{array}{lclclcl}
\mbox{\textbf{vector field}}&&\mbox{\textbf{invariant
function}}&&\mbox{\textbf{invariant
transformations}}&&\mbox{\textbf{similarity transformations}}\\
\textbf{v}_1&&\Phi(y,U,V,P,T)&&x=r+\varepsilon,y=s&&x=r,y=s,\\
\textbf{v}_2&&\Phi(x,U,V,P,T)&&x=r,y=s+\varepsilon&&x=r,y=s\\
\textbf{v}_3&&\Phi(x,y,U,V,P)&&T=f(r,s)+\varepsilon&&\mbox{translating
of all function $T$
is the similarity}\\
&&&&&&\mbox{transformation}\\
\textbf{v}_4&&\Phi(x,y,U,V,P)&&T=e^\varepsilon
f(r,s)&&\mbox{scaling of all function $T$
is the similarity}\\
&&&&&&\mbox{transformation}\\
\end{array}
\end{eqnarray*}
In the above table, the invariant solution with respect to
$\textbf{v}_1$ and $\textbf{v}_2$ obtained from translation on
both dependent variables, and for $\textbf{v}_3$ and
$\textbf{v}_4$ all invariant solutions are obtained by translating
and scaling on function $T$.
\section{Optimal system of steady two-dimensional boundary-layer stagnation-point flow equations}
Let a system of differential equation $\Delta$ admitting the
symmetry Lie group $G$,be given. Now $G$ operates on the set of
solutions $S$ of $\Delta$. Let $s\cdot G$ be the orbit of $s$, and
$H$ be a subgroup of $G$. Invariant $H-$solutions $s\in S$ are
characterized by equality $s\cdot S=\{s\}$. If $h\in G$ is a
transformation and $s\in S$,then $h\cdot(s\cdot H)=(h\cdot s)\cdot
(hHh^{-1})$. Consequently,every invariant $H-$solution $s$
transforms into an invariant $hHh^{-1}-$solution (Proposition 3.6
of \cite{[5]}).

Therefore, different invariant solutions are found from similar
subgroups of $G$. Thus, classification of invariant $H-$solutions
is reduced to the problem of classification of subgroups of $G$,up
to similarity. An optimal system of $s-$dimensional subgroups of
$G$ is a list of conjugacy inequivalent $s-$dimensional subgroups
of $G$ with the property that any other subgroup is conjugate to
precisely one subgroup in the list. Similarly, a list of
$s-$dimensional subalgebras forms an optimal system if every
$s-$dimensional subalgebra of $\goth g$ is equivalent to a unique
member of the list under some element of the adjoint
representation: $\tilde{\goth h}={\rm Ad}(g)\cdot{\goth h}$.

Let $H$ and $\tilde{H}$ be connected, $s-$dimensional Lie
subgroups of the Lie group $G$ with corresponding Lie subalgebras
${\goth h}$ and $\tilde{\goth h}$ of the Lie algebra ${\goth g}$
of $G$. Then $\tilde{H}=gHg^{-1}$ are conjugate subgroups if and
only $\tilde{\goth h}={\rm Ad}(g)\cdot{\goth h}$ are conjugate
subalgebras (Proposition 3.7 of \cite{[5]}). Thus,the problem of
finding an optimal system of subgroups is equivalent to that of
finding an optimal system of subalgebras, and so we concentrate on
it.
\subsection{One-dimensional optimal system}
For one-dimensional subalgebras,the classification problem is
essentially the same as the problem of classifying the orbits of
the adjoint representation, since each one-dimensional subalgebra
is determined by a nonzero vector in Lie algebra symmetries of
steady two-dimensional boundary-layer stagnation-point flow
equations and so to "simplify" it as much as possible.

The adjoint action is given by the Lie series
\begin{eqnarray}\label{eq:9}
\mbox{Ad}(\exp(\varepsilon\textbf{v}_i)\textbf{v}_j)=\textbf{v}_j-\varepsilon[\textbf{v}_i,\textbf{v}_j]
+\frac{\varepsilon^2}{2}[\textbf{v}_i,[\textbf{v}_i,\textbf{v}_j]]-\cdots,
\end{eqnarray}
where $[\textbf{v}_i,\textbf{v}_j]$ is the commutator for the Lie
algebra, $\varepsilon$ is a parameter, and $i,j=1,\cdots,4$. Let
$F^{\varepsilon}_i:{\goth g}\rightarrow{\goth g}$ defined by
$\textbf{v}\mapsto\mbox{Ad}(\exp(\varepsilon\textbf{v}_i)\textbf{v})$
is a linear map, for $i=1,\cdots,4$. The matrices
$M^\varepsilon_i$ of $F^\varepsilon_i$, $i=1,\cdots,4$, with
respect to basis $\{\textbf{v}_1,\cdots,\textbf{v}_4\}$ are
\begin{eqnarray}\label{eq:6}
\begin{array}{lclclcl}
M^\varepsilon_1=M^\varepsilon_2=\small\left(\begin{array}{cccc}
1&0&0&0\\0&1&0&0\\0&0&1&0\\0&0&0&1\end{array} \right),&&
M^\varepsilon_3=\small\left(\begin{array}{cccc}
1&0&0&0\\0&1&0&0\\0&0&1&0\\0&0&-\varepsilon&1\end{array}
\right),&& M^\varepsilon_4=\small\left(\begin{array}{cccc}
1&0&0&0\\0&1&0&0\\0&0&e^\varepsilon&0\\0&0&0&1\end{array} \right),
\end{array}
\end{eqnarray}
by acting these matrices on a vector field $\textbf{v}$
alternatively we can  show that a one-dimensional optimal system
of ${\goth g}$ is given by
\begin{eqnarray}\label{eq:3}
\begin{array}{lcl}
X_1=\alpha_1\textbf{v}_1+\alpha_2\textbf{v}_2,&&X_2=\alpha_1\textbf{v}_1+\alpha_2\textbf{v}_2+\alpha_3\textbf{v}_3,
\end{array}
\end{eqnarray}
where $\alpha_i$'s are real constants.
\subsection{Two-dimensional optimal system}
Next step is to construct two-dimensional optimal system, i.e.,
classification of two-dimensional subalgebras of $\goth g$. The
process is by selecting one of the vectors in (\ref{eq:3}), say,
any vector of (\ref{eq:3}). Let us consider $X_1$ (or $X_2$).
Corresponding to it, a vector field
$X=a_1\textbf{v}_1+\cdots+a_4\textbf{v}_4$, where $a_i$'s are
smooth functions of $(x,y,U,V,P,T)$ is chosen, so we must have
\begin{eqnarray}\label{eq:4}
[X_1,X]=\lambda X_1+\mu X,
\end{eqnarray}
the equation (\ref{eq:4}) leads us to the system
\begin{eqnarray}\label{eq:5}
C^i_{jk}\alpha_ja_k=\lambda
a_i+\mu\alpha_i\hspace{2cm}(i=1,2,3,4).
\end{eqnarray}
The solutions of the system (\ref{eq:5}), give one of the
two-dimensional generator and the second generator is $X_1$ or,
$X_2$ if selected. After the construction of all two-dimensional
subalgebras, for every vector fields of (\ref{eq:3}), they need to
be simplified by the action of (\ref{eq:6}) in the manner
analogous to the way of one-dimensional optimal system.

Consequently the two-dimensional optimal system of $\goth g$ has
three classes of $\goth g$'s members combinations such as
\begin{align}\label{eq:7}
&1)\;\;\beta_1\textbf{v}_1+\beta_2\textbf{v}_2,\beta_3\textbf{v}_1+\beta_4\textbf{v}_2+\textbf{v}_3,&& \beta_1\beta_3\neq\beta_2\beta_4,\\
&2)\;\;\beta_1\textbf{v}_1+\beta_2\textbf{v}_2,\textbf{v}_4,&& \beta_1^2+\beta_2^2\neq0,\nonumber\\
&3)\;\;\textbf{v}_3,\textbf{v}_4.\hfill\ \mbox{ }\nonumber
\end{align}
\subsection{Three-dimensional optimal system}
This system can be developed by the method of expansion of
two-dimensional optimal system. For this take any two-dimensional
subalgebras of (\ref{eq:7}), let us consider the first two vector
fields of (\ref{eq:7}), and rename them $Y_1$ and $Y_2$, thus, we
have a subalgebra with basis $\{Y_1,Y_2\}$, find a vector field
$Y=a_1\textbf{v}_1+\cdots+a_4\textbf{v}_4$, where $a_i$'s are
smooth functions of $(x,y,U,V,P,T)$, such the triple
$\{Y_1,Y_2,Y\}$ generates a basis of a three-dimensional algebra.
For that it is necessary an sufficient that the vector field $Y$
satisfies the equations
\begin{eqnarray}\label{eq:8}
[Y_1,Y]&=&\lambda_1Y+\mu_1Y_1+\nu_1Y_2,\nonumber\\
{[Y_2,Y]}&=&\lambda_2Y+\mu_2Y_1+\nu_2Y_2,
\end{eqnarray}
and following from (\ref{eq:8}), we obtain the system
\begin{eqnarray}\label{eq:10}
C^i_{jk}\beta_r^ja_k&=&\lambda_1a_i+\mu_1\beta_r^i+\nu_1\beta_s^i,\hspace{1cm}r=1,2,\;s=3,4,\;\alpha=1,2,3,4,\\
C^i_{jk}\beta_s^ja_k&=&\lambda_2a_i+\mu_2\beta_r^i+\nu_2\beta_s^i,\hspace{1cm}r=1,2,\;s=3,4,\;\alpha=1,2,3,4.\nonumber
\end{eqnarray}
The solutions of system (\ref{eq:10}) is linearly independent of
$\{Y_1,Y_2\}$ and give a three-dimensional subalgebra. This
process is used for the second two couple vector fields of
(\ref{eq:10}).

Consequently the three-dimensional optimal system of $\goth g$ is
given by
\begin{align*}
&1)\;\;\textbf{v}_1,\;\textbf{v}_2,\;\textbf{v}_3, &&
2)\;\;\textbf{v}_1,\;\textbf{v}_2,\;\textbf{v}_4, \\
&3)\;\;\textbf{v}_1,\;\textbf{v}_3,\;\textbf{v}_4, &&
4)\;\;\textbf{v}_2,\;\textbf{v}_3,\;\textbf{v}_4.
\end{align*}

All previous calculations lead to the table (\ref{table:2}) for
the optimal system of $\goth g$.
\begin{table}
\caption{Commutation relations of $\goth g$ }\label{table:2}
\vspace{-0.3cm}\begin{eqnarray*}\hspace{-0.75cm}\begin{array}{lclclclclc}
\hline
  dimension       &\hspace{1.1cm}1                                                                             &\hspace{2.5cm}2                                                                                                                &\hspace{0.5cm}3                                                          &\hspace{1.5cm}4  \\ \hline
                  &\hspace{1.1cm} \langle\alpha_1\textbf{v}_1+\alpha_2\textbf{v}_2\rangle                      &\hspace{0.5cm} \langle\beta_1\textbf{v}_1+\beta_2\textbf{v}_2,\beta_3\textbf{v}_1+\beta_4\textbf{v}_2+\textbf{v}_3,\rangle     &\hspace{0.5cm}\langle\textbf{v}_1,\textbf{v}_2,\textbf{v}_3\rangle       &\hspace{0.5cm}\langle\textbf{v}_1,\textbf{v}_2,\textbf{v}_3,\textbf{v}_4\rangle  \\
  subalgebras     &\hspace{1.1cm} \langle\alpha_1\textbf{v}_1+\alpha_2\textbf{v}_2+\alpha_3\textbf{v}_3\rangle &\hspace{0.5cm} \langle\beta_1\textbf{v}_1+\beta_2\textbf{v}_2,\textbf{v}_4\rangle                                              &\hspace{0.5cm}\langle\textbf{v}_1,\textbf{v}_2,\textbf{v}_4\rangle       &\hspace{0.5cm}  \\
                  &\hspace{1.1cm}                                                                              &\hspace{0.5cm} \langle\textbf{v}_3,\textbf{v}_4\rangle                                                                         &\hspace{0.5cm}\langle\textbf{v}_1,\textbf{v}_3,\textbf{v}_3\rangle       &\hspace{0.5cm}  \\
                  &\hspace{1.1cm}                                                                              &\hspace{0.5cm}                                                                                                                 &\hspace{0.5cm}\langle\textbf{v}_2,\textbf{v}_3,\textbf{v}_4\rangle       &\hspace{0.5cm}  \\
  \hline\end{array}\end{eqnarray*}\end{table}
\section{Lie Algebra Structure}
$\goth g$ has a no any non-trivial \textit{Levi decomposition} in
the form of ${\goth g}={\goth r}\ltimes{\goth g}_1$, because
$\goth g$ has no any non-trivial radical, i.e., if $\goth r$ be
the radical of $\goth g$, then ${\goth g}={\goth r}$.

If we want to integration an involuting distribution, the process
decomposes into two steps:
\begin{itemize}
\item integration of the evolutive distribution with symmetry Lie
algebra ${\goth g}/{\goth r}$, and
\item integration on integral manifolds with symmetry algebra $\goth
r$.
\end{itemize}
First, applying this procedure to the radical $\goth r$ we
decompose the integration problem into two parts: the integration
of the distribution with semisimple algebra ${\goth g}/{\goth r}$,
then the integration of the restriction of distribution to the
integral manifold with the solvable symmetry algebra $\goth r$.\\

The last step can be performed by quadratures. Moreover, every
semisimple Lie algebra ${\goth g}/{\goth r}$ is a direct sum of
simple ones which are ideal in ${\goth g}/{\goth r}$. Thus, the
Lie-Bianchi theorem reduces the integration problem ti evolutive
distributions equipped with simple algebras of symmetries. Thus,
integrating of system (\ref{eq:1}), become so much easy.

Both $\goth g$ is a solvable, because if ${\goth
g}^{(1)}=\langle\textbf{v}_i,[\textbf{v}_i,\textbf{v}_j]\rangle=[\goth
g, \goth g]$, we have ${\goth g}^{(1)}=[{\goth g},{\goth g}]
=\langle \textbf{v}_1,\cdots,\textbf{v}_4\rangle$, and ${\goth
g}^{(2)}=[{\goth g}^{(1)},{\goth g}^{(1)}] =\langle
\textbf{v}_3\rangle$, so, we have a chain of ideals ${\goth
g}^{(1)}\supset{\goth g}^{(2)}\supset\{0\}$. According to the
table of the commutators we can decompose $\goth g$ in the direct
sum of ${\goth g}=\Bbb{R}^2\oplus{\goth a}(1)$, where ${\goth
a}(1)$ is the Lie algebra of Affine transformations group $A(1)$.
\section{Conclusion}
In this article group classification of steady two-dimensional
boundary-layer stagnation-point flow equations and the algebraic
structure of the symmetry group is considered. Classification of
r-dimensional subalgebra is determined by constructing
r-dimensional optimal system. Some invariant objects are fined and
the Lie algebra structure of symmetries is found.

\end{document}